\documentclass[a4paper, 12pt]{article}

\usepackage{amsfonts, amsmath, amsthm}
\usepackage{latexsym}

\theoremstyle{plain}
\newtheorem{thm}{Theorem}[section]

\newtheorem{lem}{Lemma}[section]

\theoremstyle{definition}

\theoremstyle{remark}
\newtheorem{rmk}{Remark}[section]

\numberwithin{equation}{section}

\newcommand{\N}{\mathbb{N}}
\newcommand{\Z}{\mathbb{Z}}

\newcommand{\R}{\mathbb{R}}
\newcommand{\C}{\mathbb{C}}

\newcommand{\Hyp}{\mathbb{H}_+}

\newcommand{\Fjcub}{F_j^{\text{cub}}}
\newcommand{\conj}[1]{\overline{#1}}
\newcommand{\Lin}{L}
\newcommand{\Non}{G}

\newcommand{\pa}{\partial}
\newcommand{\eps}{\varepsilon}
\newcommand{\jb}[1]{\langle #1 \rangle}

\newcommand{\dal}{\Box}

\DeclareMathOperator{\realpart}{\rm Re}
\DeclareMathOperator{\imagpart}{\rm Im}
\DeclareMathOperator{\supp}{\rm supp}

\begin{document}
\title{
 Remarks on decay of small solutions to systems of Klein-Gordon equations with 
 dissipative nonlinearities 
 }

\author{
          Donghyun Kim\thanks{
              Department of Mathematics, Graduate School of Science, 
              Osaka University. 
              1-1 Machikaneyama-cho, Toyonaka, Osaka 560-0043, Japan. 
              (E-mail: {\tt u553252d@ecs.osaka-u.ac.jp})
             }
           \and  
          Hideaki Sunagawa\thanks{
              Department of Mathematics, Graduate School of Science, 
              Osaka University. 
              1-1 Machikaneyama-cho, Toyonaka, Osaka 560-0043, Japan. 
              (E-mail: {\tt sunagawa@math.sci.osaka-u.ac.jp})
             }
}

\date{\today }   
\maketitle

\noindent{\bf Abstract:}\ We consider the Cauchy problem for 
systems of cubic nonlinear Klein-Gordon equations in one space dimension. 
Under a  suitable structural condition on the nonlinearity, we will show that 
the small amplitude solution gains an additional logarithmic decay in 
comparison with the free evolution in the sense of $L^p$, $2\le p \le \infty$. 
\\

\\

\\


\section{Introduction}  \label{sec_intro}
We consider the Cauchy problem for
\begin{align}
 (\Box + 1)u_j = F_j(u,\pa_t u, \pa_x u), 
 \qquad (t,x) \in (0,\infty) \times \R,\ j=1,\ldots, N,
 \label{eq_NLKG}
\end{align}
with  
\begin{align}
 u_j(0,x) = \eps f_j(x), \;\; \pa_t u_j (0,x) = \eps g_j(x), 
 \qquad x \in \R,\ j=1,\ldots, N,
 \label{eq_data}
\end{align}
where $\Box = \pa_t^2 - \pa_x^2$, 
$u=(u_j)_{1\le j \le N}$ is an $\R^N$-valued unknown fuction of 
$(t,x) \in [0,\infty)\times \R$,  and the nonlinear term 
$F=(F_j)_{1\le j \le N}:\R^{3N}\to \R^N$ is assumed to be smooth and cubic 
around the origin, i.e., 
$$
 F_j(u,\pa_t u, \pa_x u) = O((|u|+|\pa_t u|+|\pa_x u|)^3)
 \quad \mbox{as}\ \ (u,\pa_t u, \pa_x u)\to (0,0,0).
$$ 
For simplicity, we always suppose $f_j$, $g_j \in C_0^{\infty}(\R)$ and 
$\eps >0$ is small enough.  We are interested in large-time behavior of the 
solution to \eqref{eq_NLKG}--\eqref{eq_data}.

From the perturbative point of view, cubic nonlinear Klein-Gordon equations in 
one space dimension are of special interest because the large-time behavior 
of the solution is actually affected by the structure of the nonlinearities 
even if the data are sufficiently small, smooth and localized. 
Let us recall some previous results briefly. 
First we focus on the single case ($N=1$). In the single case, 
a small data blow-up result is obtained by Yordanov~\cite{Yo} 
with a suitable choice of the cubic nonlinearity $F$ (see also 
Proposition~{7.8.8} in \cite{Hor}).  
On the  other hand, some classes of the cubic terms are found 
by Moriyama~\cite{Mo} and Katayama~\cite{Ka} for which the Cauchy problem 
admits a unique global solution with a free  profile if the data are small, 
smooth and compactly-supported. Their results are extended by 
Delort~\cite{De}. More precisely, a sufficient condition for small data  
global existence is introduced in \cite{De}, and the large-time asymptotic 
profile of the global solution is also provided. In the simplest case $F=u^3$, 
the result in \cite{De} can be read as follows: For  the  solution 
to $(\Box +1)u=u^3$ with $C_0^{\infty}$-data  of the size $\eps$, 
it holds that
\begin{align} 
 u(t,x)
 =
 \frac{1}{\sqrt{t}} \realpart\left(
 a(x/t) 
 e^{i\{(t^2-x^2)_+^{1/2} -\frac{3}{8} (1-|x/t|^2)_+^{1/2} |a(x/t)|^2 \log t\}}
 \right)
 +o(t^{-1/2})
 \label{eq_phase_correction}
\end{align}
as $t \to \infty$, uniformly in $x \in \R$, where $i=\sqrt{-1}$, 
$(\, \cdot\, )_+=\max\{\cdot, 0\}$, and $a(y)$ is a suitable 
$\C$-valued smooth function of $y=x/t$ vanishing when $|y|\ge 1$.  
An important consequence of this asymptotic expression is that 
the solution decays like $O(t^{-1/2})$ in $L^{\infty}$ while it does not 
behave like the  free solution. (Indeed, the solution to $(\Box+1)u_0=0$ 
with $C_0^{\infty}$-data behaves like
$$
 u_0(t,x)
 =
 \frac{1}{\sqrt{t}} \realpart\left(
 a_0(x/t)  e^{i (t^2-x^2)_+^{1/2} } \right)
 +o(t^{-1/2}) 
\qquad  \mbox{as}\ \ t \to \infty
$$
with a suitable $a_0(y)$, which can be obtained by using the stationary phase 
method; see e.g., Section~{7.2} in \cite{Hor}.) In other words, 
the additional logarithmic factor of the phase in \eqref{eq_phase_correction} 
reflects the long-range character of the cubic nonlinearity in one-dimensional 
NLKG. 
Another important class of cubic nonlinearities is the nonlinear dissipation, 
whose typical example is $F=-(\pa_t u)^3$. 
This case  is  not covered by the result in \cite{De}. 
However, as pointed out in~\cite{Su4}, one can modify the argument of 
\cite{De} to obtain the asymptotic profile also in this case. 
Moreover, a structural condition of the  cubic nonlinearity $F$ is 
introduced in \cite{Su4} under which the small data solution $u(t)$ decays 
like $O(t^{-(1/2-1/p)}(\log t)^{-1/2})$ as $t  \to \infty$ in $L^p$, 
$2\le p  \le \infty$. Remember that the $L^p$-decay rate of the  free 
evolution is $t^{-(1/2-1/p)}$. Therefore this  gain of additional logarithmic 
time decay should be interpreted as another kind of the long-range effect. 
For more recent works and related information on the single cubic nonlinear 
Klein-Gordon equations in one space dimension, we refer the readers  to  
\cite{HN1}, \cite{HN2}, \cite{LS1}, \cite{LS2} and the references 
cited therein. 

Next we turn our attention to the case of systems ($N\ge 2$), 
where the situation becomes much more complicated. 
If the linear part of \eqref{eq_NLKG} is replaced by $\Box +m_j^2$, 
a structural condition of the nonlinearities $F_j$ and 
the masses $m_j$ for small data global existence with a free profile are 
studied in \cite{Su1}. However, the result in \cite{Su1} does not 
cover the case where $m_1=\cdots=m_N$, i.e., \eqref{eq_NLKG}. 
It should be also noted that the decay rate $O(t^{-1/2})$ in $L^{\infty}$ is 
not a trivial one when $N\ge 2$; 
for instance, it is trivial to see that  the two-component system 
$$
 \left\{
 \begin{array}{l}
  (\Box+1)u_1=0, \\
  (\Box +1)u_2=u_1^3
 \end{array}
 \right.
$$
admits a global solution, but according to \cite{Su2}, the second component 
$u_2$ decays no faster than $t^{-1/2}\log t$ in $L^{\infty}$. 
In \cite{Su3}, the two-component system
\begin{align}
 \left\{
 \begin{array}{l}
  (\Box+1)u_1=(u_1^2+u_2^2)u_1, \\
  (\Box +1)u_2=(u_1^2+u_2^2)u_2
 \end{array}
 \right.
 \label{eq_NLKG2}
\end{align}
is investigated in detail, and the sharp $L^{\infty}$-decay of order 
$O(t^{-1/2})$ is obtained, but it is still the open problem to specify 
the pointwise asymptotic profile for the solution of this system. 
Note that the system \eqref{eq_NLKG2} can  be rewritten as the 
complex-valued single equation 
$(\Box +1)U=|U|^2 U$ 
through the standard identification $\R^2 \simeq \C$ (i.e., viewing 
$u=(u_1,u_2)$ as $U=u_1+i u_2$), 
however, it seems difficult to obtain the asymptotic profile for 
\eqref{eq_NLKG2} by a variant of the method of \cite{De}. 
The situation is the same in the case of nonlinear dissipation, 
that is, the previous approach \cite{Su4} does  not work well in 
the complex-valued case. A typical example of nonlinear dissipation is 
\begin{align}
 (\Box +1)U=\mu_1 |U|^2U -\mu_2 |\pa_t U|^2 \pa_t U
 \label{eq_NLKG3}
\end{align}
with $\mu_1 \in \R$ and $\mu_2>0$.
To the authors' knowledge, 
there are no previous papers which assert that the solution to 
\eqref{eq_NLKG3}
decays strictly faster than $t^{-1/2}$ in the sense of $L^{\infty}$ 
in the complex-valued case.

The aim of this paper is to introduce a structural condition of the cubic 
nonlinearity under which the solution to \eqref{eq_NLKG}--\eqref{eq_data} 
admits a global solution and it decays like $O(t^{-(1/2-1/p)}(\log t)^{-1/2})$ 
in $L^p$, $2\le p \le \infty$ if $\eps$ is small enough. As an application, 
we will see that the decay of the complex-valued solution to 
\eqref{eq_NLKG3} is strictly faster than the free evolution if $\mu_2>0$.

\section{Main Results} \label{sec_main}
In order to state the results, let us introduce some notations. 
For $j=1, \ldots, N$, we denote by $F_j^{\text{cub}}$ the cubic homogeneous 
part of $F_j$, that is, 
$$ 
 \Fjcub (\xi , \eta, \zeta)
 =
 \lim_{r \downarrow 0} r^{-3} F_j (r\xi,r\eta,r\zeta)
$$
for $(\xi, \eta, \zeta) \in \R^{3N}$. 
Roughly saying, $\Fjcub (u, \pa_t u, \pa_x u)$ gives the main part of the 
nonlinearity while $F_j (u, \pa_t u, \pa_x u) - \Fjcub (u, \pa_t u, \pa_x u)$ 
is regarded as a quartic or higher order remainder, if we are interested in 
small amplitude solutions.
Next, we set the upper branch of the unit hyperbola
$$
 \Hyp
 = 
 \{\omega = (\omega_0, \omega_1) \in \R^2 \; : \; 
    \omega_0 >0, \; \omega_0^2 - \omega_1^2 =1 \}
$$
and define $\Phi = (\Phi_j)_{1\le j\le N} : \C^N \times \Hyp \to \C^N$ by
$$
 \Phi_j (Y, \omega) 
 = 
 \frac{1}{2\pi i} \oint_{|\lambda|=1} 
 \Fjcub 
 \left(
  \realpart (Y\lambda), \; 
  -\omega_0 \imagpart (Y\lambda), \, 
  \omega_1 \imagpart(Y\lambda)
 \right)\,
 \frac{d\lambda}{\lambda^2}
$$
for $Y = (Y_j)_{1\le j\le N} \in \C^N$ and $\omega \in \Hyp$, where
$\realpart (Y\lambda)$ $=$ $\left( \realpart(Y_j \lambda) \right)_{1\le j\le N}$ and
$\imagpart (Y\lambda)$ $=$ $\left( \imagpart(Y_j \lambda) \right)_{1\le j\le N}$.
Also we denote by $\jb{\cdot, \cdot}_{\C^N}$ 
the standard scalar product in $\C^N$, i.e.,
$$
 \jb{Y, Z}_{\C^N} 
 = 
 \sum_{j=1}^N Y_j \conj{Z_j}
$$
for $Y=(Y_j)_{1\le j\le N}$, $Z=(Z_j)_{1\le j\le N} \in \C^N$. 
We also write $|Y| = \sqrt{\jb{Y, Y}_{\C^N}}$ as usual.\\\\
Now we state our main results.\\

\begin{thm}\label{thm_SDGE}
 Assume there exists a $N \times N$ positive Hermitian matrix $A$
 such that
 \begin{eqnarray} \label{eq_decayassumption0}
  \imagpart \jb{\Phi(Y, \omega), \; AY }_{\C^N} 
  \le 0
 \end{eqnarray}
 for all $(Y, \omega) \in \C^N \times \Hyp$.
 Then there exists $\eps_0 >0$ such that 
 if $\eps \in (0, \eps_0]$,
 the Cauchy problem \eqref{eq_NLKG}--\eqref{eq_data} admits 
 a unique global classical solution.
 Moreover, it satisfies
 \begin{eqnarray} \label{eq_decay0}
  \sum_{|I|\le 1}\|\pa_{t,x}^I u(t, \cdot)\|_{L^p} 
  \le 
  C(1+t)^{-(1/2 -1/p)}
 \end{eqnarray}
 for all $t \ge 0$ and $2\le p \le \infty$.\\
\end{thm}

\begin{thm}\label{thm_decay}
 Assume there exist a $N \times N$ positive Hermitian matrix $A$ 
 and a constant $C^*>0$ such that
 \begin{eqnarray} \label{eq_decayassumption1}
  \imagpart \jb{\Phi(Y, \omega), \; AY}_{\C^N} 
  \le - C^* \omega_0 |Y|^4
 \end{eqnarray}
 for all $(Y, \omega) \in \C^N \times \Hyp$. 
 Then the global solution of
 \eqref{eq_NLKG}--\eqref{eq_data} satisfies
 \begin{eqnarray} \label{eq_decay1}
  \|u(t, \cdot)\|_{L^p} 
  \le 
  C \frac{(1+t)^{-(1/2 -1/p)}}{\sqrt{\log (2+t)}}
 \end{eqnarray}
 for all $t \ge 0$ and $2\le p \le \infty$.\\
\end{thm}

\begin{thm}\label{thm_decay2}
 Assume there exist a $N \times N$ positive Hermitian matrix $A$ 
 and a constant $C_*>0$ such that
 \begin{eqnarray} \label{eq_decayassumption2}
  \imagpart \jb{\Phi(Y, \omega), \; AY}_{\C^N} 
  \le - C_* \omega_0^3 |Y|^4
 \end{eqnarray}
 for all $(Y, \omega) \in \C^N \times \Hyp$. 
 Then the global solution of
 \eqref{eq_NLKG}--\eqref{eq_data} satisfies
 \begin{eqnarray} \label{eq_decay2}
  \sum_{|I|\le 1}\|\pa_{t,x}^I u(t, \cdot)\|_{L^p} 
  \le 
  C \frac{(1+t)^{-(1/2 -1/p)}}{\sqrt{\log (2+t)}}
 \end{eqnarray}
 for all $t \ge 0$ and $2\le p \le \infty$.\\
\end{thm}

\noindent
Now we give an application of Theorem~\ref{thm_decay2} to the equation \eqref{eq_NLKG3}.
The equivalent real two-component system is 
\begin{align*}
 \left\{
 \begin{array}{l}
  (\Box+1)u_1
  =
  \mu_1(u_1^2+u_2^2)u_1 - \mu_2\{(\pa_t u_1)^2 + (\pa_t u_2)^2\} \pa_t u_1, 
  \\
  (\Box +1)u_2
  =
  \mu_1(u_1^2+u_2^2)u_2 - \mu_2\{(\pa_t u_1)^2 + (\pa_t u_2)^2\} \pa_t u_2.
 \end{array}
 \right.
\end{align*}
For this system, we have 
 \begin{eqnarray*}
 \Phi_1 (Y, \omega) = \frac{\mu_1 - i \mu_2 \omega_0^3}{8} 
 \Big{(} 3|Y_1|^2 Y_1 + 2 |Y_2|^2 Y_1 + Y_2^2 \conj{Y_1} \Big{)}, \\
 \Phi_2 (Y, \omega) = \frac{\mu_1 - i \mu_2 \omega_0^3}{8} 
 \Big{(} 3|Y_2|^2 Y_2 + 2 |Y_1|^2 Y_2 + Y_1^2 \conj{Y_2} \Big{)},
 \end{eqnarray*}
whence 
$$
 \imagpart \jb{\Phi(Y, \omega), \; AY}_{\C^2} 
 =
 -\frac{\mu_2 \omega_0^3}{8}
 \left( 2|Y_1|^4 + 2|Y_2|^4 + 4|Y_1|^2 |Y_2|^2 + |Y_1^2 + Y_2^2|^2 \right)
$$
with $A=\text{diag}(1,1)$, the identity matrix. 
This implies the assumption \eqref{eq_decayassumption2} 
of Theorem~\ref{thm_decay2} is satisfied if $\mu_2>0$. 
Therefore we conclude that 
$$
 \sum_{|I|\le 1}\|\pa_{t,x}^I U(t, \cdot)\|_{L^p} 
 \le 
 C \frac{(1+t)^{-(1/2 -1/p)}}{\sqrt{\log (2+t)}}
$$
for $2\le p \le \infty$, as desired. 
\\

\begin{rmk}\label{rmk_remark1}
An example of the system satisfying \eqref{eq_decayassumption1}
but violating \eqref{eq_decayassumption2} is
\begin{align*}
 \left\{
 \begin{array}{l}
  (\Box+1)u_1
  =
  -(u_1^2 + u_2^2) \pa_t u_1, 
  \\
  (\Box +1)u_2
  =
  -(u_1^2 + u_2^2) \pa_t u_2
 \end{array}
 \right.
\end{align*}
(or equivalently, $(\Box+1)U=-|U|^2 \pa_t U$).
For this system, we have 
 \begin{eqnarray*}
 \Phi_1 (Y, \omega) = 
  - \frac{i\omega_0}{8}
 \Big{(} |Y_1|^2 Y_1 + 2|Y_2|^2 Y_1 - Y_2^2 \conj{Y_1} \Big{)}, \\
 \Phi_2 (Y, \omega) = 
  - \frac{i\omega_0}{8}
 \Big{(} |Y_2|^2 Y_2 + 2|Y_1|^2 Y_2 - Y_1^2 \conj{Y_2} \Big{)},
 \end{eqnarray*}
so that
$$
 \imagpart \jb{\Phi(Y, \omega), \; AY}_{\C^2} 
 =
 - \frac{\omega_0}{8}
 \left( 4|Y_1|^2 |Y_2|^2 + |Y_1^2 - Y_2^2|^2 \right)
$$
with $A = \text{diag}(1,1)$. 
Therefore we can apply Theorem~\ref{thm_decay} to see that
$u(t)$ decays like $O \left( t^{-(1/2-1/p)}(\log t)^{-1/2} \right)$ as 
$t \to \infty$ in the sense of $L^p$, $2\le p\le \infty$.
However, the above theorems do not answer the question whether or not 
the same estimates hold for $\pa_{t,x} u(t, \cdot)$.\\
\end{rmk}

The rest of this paper is organized as follows. 
In Section~\ref{sec_reduction}, we reduce the original problem 
\eqref{eq_NLKG}--\eqref{eq_data} by using hyperbolic coordinates.  
Section~\ref{sec_apriori} is devoted to getting a suitable a priori estimate, 
from which the small data global existence in Theorem~\ref{thm_SDGE} 
immediately. 
After that, we prove the time decay estimates 
\eqref{eq_decay0}, \eqref{eq_decay1} and \eqref{eq_decay2} 
in Section~\ref{sec_proof}. 
In what follows, all non-negative constants will be denoted by $C$ unless 
otherwise specified.

\section{Reduction of the Problem} \label{sec_reduction}
In this section, we perform some reduction of the problem along the idea of \cite{De}
with a slight modification. 
In the following, as mentioned before, we shall neglect the higher order terms of $F_j$ 
(i.e. we assume $F_j = \Fjcub$) because the higher order terms do not have an essential
influence on our problem.

Let $B$ be a positive constant which satisfies 
$$
 \supp f_j \cup \supp g_j \subset \{ x \in \R \,:\, |x| \le B \} 
$$
and let $\tau_0 > 1 + 2B$. We start with the fact that we may treat the problem 
as if the Cauchy data are given on the upper branch of the hyperbola
$$
 \{ (t,x) \in \R^2 \,:\, (t+2B)^2 - x^2 = \tau_0^2, \; t >0 \}
$$
and they are sufficiently smooth, small, compactly-supported. 
This is a consequence of the classical local existence theorem 
and the finite speed of propagation 
(see e.g., Proposition 1.4 of \cite{De} for the detail). 
Next, let us introduce the hyperbolic coordinates
$(\tau, z ) \in [\tau_0 , \infty) \times \R$ in the interior of the light cone, i.e.,
$$
 t+2B = \tau \omega_0 (z),\quad x = \tau \omega_1 (z)
$$ 
for $|x| < t + 2B$ where
$\omega(z) = (\omega_0(z), \omega_1(z))=(\cosh z, \sinh z) \in \Hyp$. 
Note that $\omega_k (z)$ satisfies
$$
 |\omega_k (z)| \le Ce^{|z|}
$$
for all $z \in \R$ and $k=0,1$.
Then we can easily check that
$$
 \pa_t = \omega_0(z) \pa_\tau - \frac{1}{\tau}\omega_1(z) \pa_z,\quad 
 \pa_x = -\omega_1(z) \pa_\tau + \frac{1}{\tau}\omega_0(z) \pa_z,\quad
 \dal = \pa_\tau^2 + \frac{1}{\tau} \pa_\tau - \frac{1}{\tau^2} \pa_z^2.
$$
We also take a weight function $\chi(z) \in C^\infty(\R)$ satisfying
$$
 0 < \chi (z) \le C_0 e^{-\kappa |z|} \;\;\text{and}\;\; 
 | \pa_z^j \chi (z) | \le C_j \chi (z) \;\;(j=1,2,\cdots)
$$
with a large parameter $\kappa \gg 1$ and positive constants $C_0, C_j$.
With this weight function, 
let us define the new unknown function $v(\tau, z) = (v_j(\tau,z))_{1\le j\le N}$ by
\begin{equation*}
 u_j (t,x)
 =
 \frac{\chi(z)}{\sqrt \tau} v_j (\tau,z).
\end{equation*}
Then we see that $v_j$ satisfies
\begin{equation*}
 \Lin v_j
 =
 \Non_j (\tau, z, v, \pa_\tau v, \pa_z v)
\end{equation*}
if $u$ solves \eqref{eq_NLKG}, where
$$
 \Lin 
 = 
 \pa_\tau^2 - \frac{1}{\tau^2} \left( 
 \pa_z^2 + 2 \frac{\chi' (z)}{\chi(z)} \pa_z + \frac{\chi''(z)}{\chi(z)} - \frac{1}{4} 
 \right) +1
$$
and
$$
 \Non_j 
 = 
 \frac{\chi(z)^2}{\tau} \Fjcub(v,\, \omega_0(z)\pa_\tau v,\, -\omega_1(z)\pa_\tau v)
 +{Q}_j(\tau, z, v, \pa_\tau v, \pa_z v).
$$
Here ${Q}_j$ takes the form of
$$
 {Q}_j
 =
 \sum_{\text{finite}} \bigg( 
 (\tau^{-2} \text{ or } \tau^{-3} \text{ or } \tau^{-4}) 
 \times (\text{bounded functions of } z) 
 \times P_3(v, \pa_\tau v, \pa_z v)  \bigg),
$$
where $P_3(v, \pa_\tau v, \pa_z v)$ is a homogeneous polynomial of order $3$ 
consisting of 
$(v_1, \cdots, v_N$, 
$\pa_\tau v_1, \cdots, \pa_\tau v_N,$ 
$\pa_z v_1, \cdots,  \pa_z v_N)$.
As we shall see in \eqref{est_Qj} below, 
$Q_j$ can be regarded as a remainder, 
while the first term of $G_j$ plays a role as a main term.

At last, the original problem \eqref{eq_NLKG}--\eqref{eq_data} is reduced to
\begin{align}
 \left\{
 \begin{array}{l}
  \Lin v_j 
  = 
  \Non_j (\tau, z, v, \pa_\tau v, \pa_z v), \quad 
  \tau > \tau_0, \; z \in \R, \; j=1, \cdots, N, 
  \\
  (v_j, \pa_\tau v_j) |_{\tau=\tau_0}
  =
  (\eps \tilde{f_j}, \eps \tilde{g_j}), \quad 
  z \in \R, \; j=1, \cdots, N 
 \end{array}
 \right.
\label{eq_reduced}
\end{align}
where $\tilde{f_j}, \tilde{g_j}$ are sufficiently smooth functions of $z$ 
with compact support.

\section{A Priori Estimate} \label{sec_apriori}
This section is devoted to getting an a priori estimate 
for the solution of the reduced problem \eqref{eq_reduced}
under the condition \eqref{eq_decayassumption0}.
First, we set
$$
 M(T)
 =
 \sup_{(\tau,z)\in[\tau_0,T)\times\R} 
 \left( |v(\tau,z)| + |\pa_\tau v(\tau,z)| + \frac{1}{\tau} |\pa_z v(\tau,z)| \right)
$$
for the smooth solution $v(\tau,z)$ to \eqref{eq_reduced} on $\tau \in [\tau_0, T)$.
We will prove the following:

\begin{lem} \label{lem_apriori}
 Under the assumption of Theorem~\ref{thm_SDGE}, 
 there exist $\eps_1 >0$ and $C_1 >0$ such that 
 $M(T) \le \sqrt \eps$ implies $M(T) \le C_1 \eps$ for any $\eps \in (0, \eps_1]$.
 Here $C_1$ is independent of $T$.
\end{lem}

Once this lemma is proved, 
we can derive the global existence part of Theorem~\ref{thm_SDGE} 
in the following way: By taking $\eps_0 \in (0, \eps_1]$ so that 
$2C_1 \eps_0^{1/2} \le 1$, we deduce that $M(T) \le \eps^{1/2}$ implies 
$M(T) \le \eps^{1/2} /2$ for any $\eps \in (0,\eps_0]$. 
Then by the continuity argument, 
we have $M(T) \le C_1 \eps$ as long as the solution exists. 
Therefore the local solution to \eqref{eq_reduced} can be extended to the global one.
Going back to the original variables, 
we deduce the small data global existence for \eqref{eq_NLKG}--\eqref{eq_data}.

The proof of Lemma~\ref{lem_apriori} will be divided into two steps: 
We first derive an auxiliary estimate for the energy 
$$
 E_s (\tau) 
 = 
 \sum_{k=0}^s \frac{1}{2} \int_{\R} 
 | \pa_\tau \pa_z^k v(\tau, z) |^2 
 + \left| \frac{\pa_z}{\tau} \pa_z^k v(\tau, z) \right|^2 
 + | \pa_z^k v(\tau, z) |^2 \, dz
$$
for $s \in \N_0$ under the assumption $M(T) \le \sqrt \eps$.
Remark that we do not need a special structure of the nonlinearity at this 
stage. 
Next we will prove the desired estimate for $M(T)$ 
by using the the structural condition \eqref{thm_SDGE} 
and the auxiliary estimate for the energy obtained in the first step. \\

\textbf{Proof of Lemma~\ref{lem_apriori}.}

\underline{(Step 1)} 
This part is essentially the same as that of the previous works. 
Our goal here is to show
\begin{eqnarray} \label{eq_energy}
 E_s (\tau)
 \le
 C\eps^2 \tau^\delta
\end{eqnarray}
under the assumption that $M(T) \le \sqrt \eps$,
where $s \ge 3$ and $0<\delta<\frac{1}{3}$.
Let us introduce
$$
 E_s (\tau; w) 
 = 
 \sum_{k=0}^s \frac{1}{2} \int_{\R} 
 | \pa_\tau \pa_z^k w(\tau, z) |^2 
 + \left| \frac{\pa_z}{\tau} \pa_z^k w(\tau, z) \right|^2 
 + | \pa_z^k w(\tau, z) |^2 \, dz
$$
for $s \in \N_0$ and for smooth function $w$ of $(\tau, z) \in [\tau_0, T) \times \R$.
We start with the following energy inequality,
whose proof is found in Appendix of \cite{Su3}
(see also \textsection 3 of \cite{Su2}, or \textsection 3 of \cite{Su4}).\\

\begin{lem} \label{lem_energy}
For $s \in \N_0$ and $l = 0,1,$ we have
$$
 \frac{d}{d \tau} E_s (\tau;w) 
 \le 
 \frac{C'}{\tau^{1+l}}E_{s+l}(\tau;w)
 +CE_s(\tau;w)^{\frac{1}{2}}\|\Lin w(\tau,\cdot)\|_{H^s}
 +\frac{C}{\tau^2}E_s(\tau;w)
$$
where $C' = 2 \sup_{z \in \R} \frac{|\chi'(z)|}{\chi(z)}$ 
and $\| \cdot \|_{H^s}$ denotes the standard norm of the Sobolev space $H^s$.\\
\end{lem}

We shall apply the above lemma with $l=0$, $s=s_0+s_1+1$ and 
$w=v_1, \cdots, v_N$,
where $s_0$ is an integer greater than $C'$ 
and $s_1$ is a fixed arbitrary non-negative integer. 
Since the Gagliardo-Nirenberg inequality yields 
$$
 \| \Non_j (\tau,\cdot,v(\tau,\cdot),\pa_\tau v(\tau,\cdot),\pa_z v(\tau,\cdot))\|_{H^s}
 \le \frac{C}{\tau} M(T)^2 E_s(\tau)^{\frac{1}{2}} 
 \le \frac{C}{\tau} \eps E_s(\tau)^{\frac{1}{2}},
$$ 
we have
$$
 \frac{d}{d\tau} E_{s_0 + s_1 +1}(\tau) 
 \le 
 \left( \frac{C' + C\eps}{\tau} + \frac{C}{\tau^2} \right) E_{s_0 + s_1 +1} (\tau) 
 \le \left(\frac{s_0+\frac{1}{2}}{\tau}+\frac{C}{\tau^2}\right)E_{s_0+s_1+1}(\tau)
$$ 
for sufficiently small $\eps$. Thus the Gronwall lemma yields 
$$
 E_{s_0+s_1+1}(\tau) 
 \le 
 E_{s_0+s_1+1}(\tau_0) 
 \exp \left(\int_{\tau_0}^\tau \frac{s_0 + \frac{1}{2}}{\sigma} 
 +\frac{C}{\sigma^2} \, d \sigma \right) \le C \eps^2 \tau^{s_0 + \frac{1}{2}}.
$$
Next, we apply Lemma~\ref{lem_energy} with $l=1, s=s_0+s_1$. 
Using the above inequality, we have 
\begin{eqnarray*}
 \frac{d}{d \tau} E_{s_0 + s_1} (\tau) 
 &\le& 
 \frac{C'}{\tau^2} E_{s_0+s_1+1}(\tau) 
 + \frac{C\eps}{\tau}E_{s_0+s_1} (\tau) 
 + \frac{C}{\tau^2}E_{s_0+s_1}(\tau) \\
 &\le& 
 C\eps^2 \tau^{s_0-\frac{3}{2}}
 +\left(\frac{C\eps}{\tau}+\frac{C}{\tau^2} \right) E_{s_0+s_1}(\tau).
\end{eqnarray*}
Therefore it follows from the Gronwall lemma that 
$$
 E_{s_0+s_1}(\tau) 
 \le 
 C\eps^2 \tau^{s_0 - \frac{1}{2}}.
$$
Repeating the same procedure $n$ times, we have 
$$
 E_{s_0+s_1+1 - n}(\tau) 
 \le C\eps^2 \tau^{s_0 - n + \frac{1}{2}}
$$
for $n=1,2,\cdots,s_0$. In particular we have 
$$
 E_{s_1 +1} (\tau) 
 \le 
 C \eps^2 \tau^{\frac{1}{2}}.
$$
Finally, we again use Lemma~\ref{lem_energy} with $l=1, s=s_1$ to obtain 
\begin{eqnarray*}
 \frac{d}{d\tau}E_{s_1}(\tau) 
 &\le& 
 \frac{C'}{\tau^2}E_{s_1+1}(\tau) 
 + \frac{C\eps}{\tau} E_{s_1}(\tau) 
 + \frac{C}{\tau^2}E_{s_1}(\tau) \\
 &\le& 
 \frac{C\eps^2}{\tau^{3/2}} 
 + \left( \frac{C\eps}{\tau} + \frac{C}{\tau^2} \right) E_{s_1}(\tau).
\end{eqnarray*}
The Gronwall lemma yields 
$$
 E_{s_1}(\tau) 
 \le C\eps^2 \exp \left( 
 \int_{\tau_0}^\tau \frac{C\eps}{\sigma}+\frac{C}{\sigma^2}\,d\sigma\right) 
 \le C\eps^2 \tau^{C\eps}
$$
for $\tau \in [\tau_0, T)$. Replacing $s_1$ by $s$ 
and choosing $\eps$ small that $C\eps \le \delta$, we arrive at 
\eqref{eq_energy}.\\

\underline{(Step 2)} Now we are going to prove $M(T) \le C_1 \eps$
under the same assumption as before.
As in \cite{Su2}, \cite{Su3} and \cite{Su4},
we introduce the $\C^N$-valued function $\alpha = (\alpha_j)_{1\le j\le N}$ by
$$
 \alpha_j (\tau, z) 
 = 
 e^{-i \tau} \left(1 + \frac{1}{i} \frac{\pa}{\pa \tau}\right) v_j (\tau, z), 
 \qquad (\tau, z) \in [\tau_0, T) \times \R, \quad j=1, \cdots, N
$$
for the solution $v(\tau,z)$ to \eqref{eq_reduced}.
In view of the relations
$$
 |\alpha(\tau,z)|
 =
 \left(|v(\tau,z)|^2 + |\pa_\tau v(\tau,z)|^2 \right)^{1/2}
$$
and
$$
 \frac{1}{\tau}|\pa_z v(\tau,z)|
 \,\le\,
 \frac{CE_2(\tau)^{1/2}}{\tau}
 \,\le\,
 C\eps \tau^{-(1-\delta/2)}
 \,\le\,
 C\eps,
$$
it suffices to show that
\begin{eqnarray} \label{eq_alphaineq}
 \sup_{(\tau,z)\in[\tau_0,T)\times\R}|\alpha(\tau,z)|
 \le
 C\eps
\end{eqnarray}
holds true under the assumptions 
$M(T)\le\sqrt \eps$ and \eqref{eq_energy}.
For this purpose,  first we note that
\begin{eqnarray} \label{eq_alphaDE1}
 \frac{\pa \alpha_j}{\pa \tau} 
 &=& 
 -ie^{-i \tau} (\pa_\tau^2 +1) v_j \nonumber \\
 &=& 
 -ie^{-i \tau} \left( \Non_j (\tau,z,v,\pa_\tau v,\pa_z v) + \frac{1}{\tau^2} \left( 
 \pa_z^2 + 2 \frac{\chi'(z)}{\chi(z)} \pa_z + \frac{\chi''(z)}{\chi(z)} - \frac{1}{4} 
 \right) v_j \right) \nonumber \\
 &=& 
 -ie^{-i \tau} \frac{\chi(z)^2}{\tau} 
 \Fjcub (v, \omega_0 (z) \pa_\tau v, - \omega_1 (z) \pa_\tau v) 
 + \frac{R_j (\tau, z)}{\tau^2},
\end{eqnarray}
where
\begin{eqnarray*}
 R_j (\tau, z) 
 &=& 
 -i e^{-i\tau} \tau^2 Q_j
 - ie^{-i\tau} \left( 
     \pa_z^2 + 2 \frac{\chi'(z)}{\chi(z)} \pa_z +\frac{\chi''(z)}{\chi(z)} 
     - \frac{1}{4} 
   \right) v_j.
\end{eqnarray*}
Using the assumption $M(T) \le \sqrt \eps$, we have 
$$
 \left| 
 \frac{\chi(z)^2}{\tau}\Fjcub(v, \omega_0(z) \pa_\tau v, -\omega_1(z) \pa_\tau v) 
 \right| 
 \,\le\, 
 \frac{C}{\tau} M(T)^3 
 \,\le\, 
 \frac{C}{\tau}\eps^{3/2}.
$$
Also we have 
\begin{eqnarray} \label{eq_remainder}
 |R_j (\tau, z) | 
 \,\le\, 
 C \eps \tau^{3\delta/2}
\end{eqnarray}
since 
\begin{eqnarray}
 |Q_j|
 \le
 \frac{C}{\tau^2}(|\pa_\tau v|^3 + |\pa_z v|^3 +|v|^3) 
 \le
 \frac{C}{\tau^2} E_2(\tau)^{3/2}  
 \le 
 C\eps^3 \tau^{3\delta/2-2}
 \label{est_Qj}
\end{eqnarray}
and 
$$
 \left| \left( 
  \pa_z^2 + \frac{2 \chi'(z)}{\chi(z)} \pa_z + \frac{\chi''(z)}{\chi(z)} - \frac{1}{4} 
 \right) v_j \right| 
 \,\le\,
 C \sum_{k=0}^2 | \pa_z^k v | 
 \,\le\,
 C E_3(\tau)^{1/2} 
 \,\le\, 
 C \eps \tau^{\delta/2}. 
$$
Combining all together, we obtain
\begin{eqnarray} \label{eq_alpha-estimate}
 \left| \frac{\pa \alpha_j}{\pa \tau} (\tau, z) \right| 
 \,\le\, 
 \frac{C\eps^{3/2}}{\tau} + \frac{C\eps}{\tau^{2-3\delta/2}} \;\le\; \frac{C\eps}{\tau}.
\end{eqnarray}
For the next step, we introduce the function $H_j(\theta,z)$ which is defined 
by 
$$
 H_j(\theta,z) 
 = 
 e^{-i \theta} \Fjcub (
 \realpart (\alpha e^{i \theta}),\, 
 -\omega_0(z) \imagpart(\alpha e^{i \theta}),\, 
 \omega_1(z) \imagpart(\alpha e^{i \theta}))
$$
with $\alpha$ being regarded as a parameter for the moment. 
Since $H_j$ is $2\pi$-periodic with respect to $\theta$, we have 
$$
 H_j(\theta,z) 
 = 
 \sum_{n \in \Z} \widehat{H}_{j, n}(z) e^{i n \theta},
$$
where $\widehat{H}_{j,n}$ denotes the $n$-th Fourier coefficient, i.e.,
$$
 \widehat{H}_{j,n}(z) 
 = 
 \frac{1}{2\pi} \int_0^{2\pi} H_j(\theta,z) e^{-i n \theta} \,d\theta.
$$
Noting that
$$
 \int_0^{2\pi} \Fjcub 
 (\realpart(\alpha e^{i \theta}),\, 
 -\omega_0(z) \imagpart(\alpha e^{i \theta}),\, 
 \omega_1(z) \imagpart(\alpha e^{i \theta})) 
 e^{-i(n+1) \theta} \,d\theta 
 =
 0
$$
when $n+1 \notin \{1,-1,3,-3\}$, we see that
$$
 H_j(\theta,z) 
 = 
 \widehat{H}_{j, 0}(z) 
 + \widehat{H}_{j,2} e^{2i \theta} 
 + \widehat{H}_{j,-2}(z) e^{-2i \theta} 
 + \widehat{H}_{j,-4} (z) e^{-4i \theta}.
$$
Also we observe that
\begin{eqnarray*}
 \widehat{H}_{j,0}(z)
 =
 \frac{1}{2\pi i} \oint_{|\lambda|=1} 
 \Fjcub 
 \left(
  \realpart (\alpha \lambda), \; 
  -\omega_0(z) \imagpart (\alpha \lambda), \, 
  \omega_1(z) \imagpart( \alpha \lambda)
 \right)\,
 \frac{d\lambda}{\lambda^2}
 =
 \Phi_j(\alpha, \omega(z)).
\end{eqnarray*}
From \eqref{eq_alphaDE1}, we have
\begin{eqnarray} \label{eq_alphaDE2}
 \frac{\pa \alpha_j}{\pa \tau}
 &=&
 -\frac{i\chi(z)^2}{\tau} e^{-i \tau} 
 \Fjcub\left(\realpart(\alpha e^{i \tau}),\,
 -\omega_0 (z) \imagpart(\alpha e^{i \tau}),\,
 \omega_1(z)  \imagpart(\alpha e^{i \tau}) \right)
 +\frac{R_j (\tau, z)}{\tau^2} \nonumber\\
 &=&
 -\frac{i \chi(z)^2}{\tau} H_j (\tau,z) + \frac{R_j (\tau, z)}{\tau^2} \nonumber\\
 &=& 
 - \frac{i \chi(z)^2}{\tau} \left(
    \widehat{H}_{j, 0}(z) 
  +\widehat{H}_{j, 2}(z) e^{ 2i \tau} 
  +\widehat{H}_{j,-2}(z) e^{-2i \tau} 
  +\widehat{H}_{j,-4}(z) e^{-4i \tau} \right)
  +\frac{R_j (\tau, z)}{\tau^2} \nonumber \\
 &=& 
 -\frac{i \chi(z)^2}{\tau} \Phi_j (\alpha, \omega(z)) 
 +S_j (\tau, z) 
 +\frac{R_j (\tau, z)}{\tau^2}, 
\end{eqnarray}
where 
$\realpart(\alpha e^{i \tau}) = (\realpart(\alpha_j e^{i\tau}))_{1\le j\le N}$,
$\imagpart(\alpha e^{i \tau}) = (\imagpart(\alpha_j e^{i\tau}))_{1\le j\le N}$
and
$$
 S_j (\tau, z) 
 = 
 \frac{\chi(z)^2}{\tau} \bigg{(}  
   I_{j,-1} (\alpha, \omega(z)) e^{-2i \tau} 
 +I_{j, 3} (\alpha, \omega(z)) e^{2i \tau} 
 +I_{j,-3} (\alpha, \omega(z)) e^{-4i \tau} \bigg{)},
$$
$$
 I_{j,n} (\alpha, \omega) 
 =
 \frac{1}{2\pi i} \int_0^{2\pi} e^{-i n \theta} 
 \Fjcub (
   \realpart(\alpha e^{i \theta}),\, 
  -\omega_0 \imagpart(\alpha e^{i \theta}),\, 
   \omega_1 \imagpart(\alpha e^{i \theta}))\,d\theta.
$$
Now we put 
$\nu_A(Y) =\sqrt{\jb{ Y, A Y}_{\C^N}}$ for $Y \in \C^N$,
where the matrix $A$ is from the assumption of Theorem~\ref{thm_SDGE}.
Then we can immediately check that
\begin{eqnarray} \label{eq_eigen}
 |\jb{Y, AZ}_{\C^N}|
 \le 
 \nu_A(Y) \nu_A(Z), 
 \qquad 
 \sqrt{\lambda_*}|Y| 
 \le 
 \nu_A(Y) 
 \le 
 \sqrt{\lambda^*}|Y|
\end{eqnarray}
for $Y, Z \in \C^N$, where $\lambda^*$ (resp. $\lambda_*$) is the 
largest (resp. smallest) eigenvalue of $A$.
With the notations 
$S = (S_j)_{1\le j\le N}$, $R = (R_j)_{1\le j\le N}$, 
it follows from 
\eqref{eq_alphaDE2},
\eqref{thm_SDGE},
\eqref{eq_eigen}
and
\eqref{eq_remainder}
that
\begin{eqnarray*}
 \frac{\pa}{\pa \tau} \left( \nu_A(\alpha(\tau,z))^2 \right) 
 &=& 
 2\realpart \jb{\pa_\tau \alpha(\tau,z), A\alpha(\tau,z)}_{\C^N} \\
 &=& 
 \frac{2\chi(z)^2}{\tau}\imagpart\jb{\Phi(\alpha,\omega(z)),A\alpha}_{\C^N} 
  +2\realpart \jb{S, A\alpha}_{\C^N} 
  +\frac{2}{\tau^2} \realpart \jb{R, A\alpha}_{\C^N} \\
 &\le& 
  2\realpart \jb{S, A \alpha}_{\C^N} 
  +\frac{2}{\tau^2}| \jb{R, A\alpha}_{\C^N} | \\
 &\le&
 2\realpart \jb{S, A \alpha}_{\C^N}
  +\frac{1}{\tau^2} \left( \nu_A(\alpha)^2 + \nu_A(R)^2 \right)\\
 &\le& 
 2\realpart \jb{S, A \alpha}_{\C^N} 
  +\frac{1}{\tau^2} \nu_A(\alpha)^2 
  +\frac{C \eps^2}{\tau^{2- 3\delta}}
\end{eqnarray*}
which yields
\begin{eqnarray}
 \nu_A(\alpha(\tau, z))^2 
 &\le& 
 C\eps^2
  +
 2 \left| \int_{\tau_0}^\tau \realpart \jb{S, A\alpha}_{\C^N}\,d\sigma \right| 
  +
 \int_{\tau_0}^\tau \nu_A(\alpha(\sigma,z))^2\,\frac{d\sigma}{\sigma^2}
 \nonumber\\
 &\le& 
 C\eps^2 
  +
 \int_{\tau_0}^\tau \nu_A(\alpha(\sigma,z))^2 \,\frac{d\sigma}{\sigma^2}
\end{eqnarray}
for $\tau \in [\tau_0, T)$, provided that 
\begin{eqnarray} \label{eq_resonant}
 \sup_{\tau \in [\tau_0, T)}  
  \left| \int_{\tau_0}^\tau \realpart \jb{S, A\alpha}_{\C^N}\,d\sigma \right| 
 \,\le\, 
 C\eps^2.
\end{eqnarray}
Once we get \eqref{eq_resonant},
we can apply the Gronwall lemma 
and \eqref{eq_eigen} to obtain \eqref{eq_alphaineq}.

It remains to prove \eqref{eq_resonant}.
To this end, we observe that
\begin{eqnarray*}
 \int_{\tau_0}^\tau
  \left(\prod_{m=1}^4 \alpha_{k_m}^{(\gamma_m)}\right)
  \frac{e^{i{b}\sigma}}{\sigma}\,d\sigma 
 &=& 
 \int_{\tau_0}^\tau (\pa_\tau K_1)(\sigma,z) + K_2(\sigma,z)\,d \sigma \\
 &=& 
 K_1 (\tau,z) - K_1 (\tau_0,z) + \int_{\tau_0}^\tau K_2 (\sigma,z)\,d\sigma,
\end{eqnarray*}
for $b \in \R\backslash\{0\}$, $k_1, \cdots, k_4 \in \{1, \cdots, N\}$ 
and $\gamma_1, \cdots, \gamma_4 \in \{+,-\}$, 
where 
$\alpha_k^{(+)}=\alpha_k$,
$\alpha_k^{(-)}=\conj{\alpha_k}$
and
$$
 K_1(\tau,z) 
 = 
 \left(\prod_{m=1}^4 \alpha_{k_m}^{(\gamma_m)}\right)
  \frac{e^{i {b} \tau}}{(i {b} \tau)}, 
$$
$$
 K_2(\tau,z) 
 = 
 \left(\prod_{m=1}^4 \alpha_{k_m}^{(\gamma_m)}\right)
  \frac{e^{i {b} \tau}}{i{b} \tau^2} 
 -
  \sum_{m=1}^4 \Bigg{(} \left(\pa_\tau \alpha_{k_m}^{(\gamma_m)}\right)
   \prod_{l \in\{1,2,3,4\} \backslash \{m\}}
   \alpha_{k_l}^{(\gamma_l)}\Bigg{)}
   \frac{e^{i {b} \tau}}{i {b}\tau}.\\ 
$$
Using \eqref{eq_alpha-estimate} and $M(T)\le\sqrt\eps$, we have 
$$
 |K_1(\tau,z)| 
 \le 
 \frac{C\eps^2}{\tau}, 
 \qquad |K_2(\tau,z)|
 \le 
 \frac{C\eps^2}{\tau^2} 
  +\frac{C\eps}{\tau} \cdot \eps^{\frac{3}{2}} \cdot \frac{1}{\tau} 
 \le 
 \frac{C\eps^2}{\tau^2}.
$$
From them we deduce that
\begin{eqnarray*}
 \sup_{\tau \in [\tau_0, T)}  
  \left| \int_{\tau_0}^\tau \realpart \jb{S, A\alpha}_{\C^N}\,d\sigma \right|
 &\le&
 \sum_{b \in \{-2,2,4\}} C
 \sup_{\tau \in [\tau_0,T)} \left| 
 \int_{\tau_0}^\tau \left(\prod_{m=1}^4 \alpha_{k_m}^{(\gamma_m)}\right)
 \frac{e^{i {b} \sigma}}{\sigma} \; d \sigma \right| \\
 &\le& 
 C\eps^2 \left( 1+ \int_1^\infty \frac{d \sigma}{\sigma^2} \right).
\end{eqnarray*}
%
%
This completes the proof of Lemma~\ref{lem_apriori}.
\qed

\section{Proof of the Decay Estimates} \label{sec_proof}
We are in a position to prove the time decay estimates 
\eqref{eq_decay0}, \eqref{eq_decay1}, \eqref{eq_decay2}
under the conditions \eqref{eq_decayassumption0},
\eqref{eq_decayassumption1}, \eqref{eq_decayassumption2}, respectively. 
First, we remember that our change of variable is
\begin{eqnarray} \label{eq_substitution}
 u_j(t,x)
 =
 \frac{\chi(z)}{\sqrt\tau} \realpart(\alpha_j(\tau,z)e^{i\tau})
\end{eqnarray}
with
$t+2B = \tau\omega_0(z)$,
$x=\tau\omega_1(z)$
for $|x| < t+2B$, 
and that $u(t,\cdot)$ is supported on $\{x \in \R :|x|\le t+B\}$.
So it follows from \eqref{eq_alphaineq} and \eqref{eq_substitution} that
\begin{eqnarray}\label{eq_udecay0}
 |u(t,x)|
 \,\le\,
 \frac{\chi(z)|\alpha(\tau,z)| \omega_0(z)^{1/2}}{\sqrt{\tau}\omega_0(z)^{1/2}}
 \,\le\,
 C\eps(1+t)^{-1/2}.
\end{eqnarray}
Using \eqref{eq_udecay0} and the finite propagation speed, we have 
\begin{eqnarray*}
 \|u(t, \cdot) \|_{L^p(\R)} 
 &=& 
 \|u(t, \cdot) \|_{L^p(\{ |x| \le t+B\})} \\
 &\le& 
 C\eps (1+t)^{-1/2}
  \left( \int_{|x| \le t+B}1\,dx \right)^{1/p}\\
 &\le& 
 C\eps (1+t)^{-(1/2-1/p)}
\end{eqnarray*}
for $p \in [2,\infty]$.
Also we note that $\pa_t u_j$ and $\pa_x u_j$ can be written as
\begin{eqnarray*}
 \pa_t u_j(t,x)
 &=&
 \left( \omega_0(z)\pa_\tau - \frac{1}{\tau}\omega_1(z)\pa_z\right)
 \left(\frac{\chi(z)}{\sqrt{\tau}}v_j(\tau,z) \right)\\
 &=&
 - \frac{\chi(z)\omega_0(z)^{3/2}}{\sqrt{\tau\omega_0(z)}}
   \imagpart(\alpha_j e^{i\tau})
 - \frac{\chi(z)\omega_0(z)}{\tau^{3/2}}
  \left( \frac{\omega_1(z)}{\omega_0(z)}\pa_z
  + \frac{\omega_1(z)}{\omega_0(z)} \frac{\chi'(z)}{\chi(z)}
  + \frac{1}{2} \right) v_j
\end{eqnarray*}
and
\begin{eqnarray*}
 \pa_x u_j(t,x)
 &=&
 \left( -\omega_1(z)\pa_\tau + \frac{1}{\tau}\omega_0(z)\pa_z\right)
 \left(\frac{\chi(z)}{\sqrt{\tau}}v_j(\tau,z) \right)\\
 &=&
 \frac{\frac{\omega_1(z)}{\omega_0(z)}\chi(z)\omega_0(z)^{3/2}}{\sqrt{\tau\omega_0(z)}}
   \imagpart(\alpha_j e^{i\tau})
 + \frac{\chi(z)\omega_0(z)}{\tau^{3/2}}
  \left( \pa_z
  +\frac{\chi'(z)}{\chi(z)}
  + \frac{1}{2} \frac{\omega_1(z)}{\omega_0(z)} \right) v_j.
\end{eqnarray*}
Thus we get
\begin{eqnarray*}
 \sum_{|I|=1} \left| \pa_{t,x}^I u(t,x) \right|
 &\le&
 \frac{C\eps}{\sqrt{\tau\omega_0(z) }}
 +\frac{C\eps\chi(z)\omega_0(z)^{1+4/3}}{(\tau\omega_0(z))^{4/3}}
 +\frac{C\eps\chi(z)\omega_0(z)^{1+3/2}}{(\tau\omega_0(z))^{3/2}}\\
 &\le& 
 C\eps (1+t)^{-1/2}
 +C\eps(1+t)^{-4/3}
 +C\eps(1+t)^{-3/2}\\
 &\le&
 C\eps (1+t)^{-1/2}
\end{eqnarray*}
and
$$
 \sum_{|I|\le 1} \| \pa_{t,x}^I u(t, \cdot) \|_{L^p(\R)}
 \le
 C \eps (1+t)^{-(1/2-1/p)}
$$
for $p \in [2,\infty]$.
This proves Theorem~\ref{thm_SDGE}.\\

Now we are going to prove Theorem~\ref{thm_decay}. 
From the argument of the previous section,
we see that
$$
 \frac{\pa \alpha_j}{\pa \tau}
 \,=\,
 -\frac{i\chi(z)^2}{\tau}\Phi_j(\alpha,\omega(z))
 +S_j(\tau,z)
 +\frac{R_j(\tau,z)}{\tau^2}
$$
and that $S_j(\tau,z)$ can be written as
(by the similar method as in the previous section)
$$
 S_j(\tau,z)
 \,=\,
 \pa_\tau V_j(\tau,z) + W_j(\tau,z)
$$
with
$$
 |V_j (\tau,z)|
 \le
 \frac{C\eps^3}{\tau},\qquad
 |W_j(\tau,z)|
 \le
 \frac{C\eps^3}{\tau^2}.
$$
Putting $\beta_j = \alpha_j - V_j$
and
$\beta = (\beta_j)_{1\le j\le N}$, 
$V = (V_j)_{1\le j\le N}$, we have
\begin{eqnarray} \label{eq_betaDE}
 \frac{\pa \beta_j}{\pa \tau} 
 \,=\, 
 -\frac{i \chi(z)^2}{\tau} \Phi_j (\beta, \omega(z))+\rho_j(\tau,z)
\end{eqnarray}
where
$$
 \rho_j(\tau,z) 
 = 
 W_j + \frac{R_j}{\tau^2} 
 -\frac{i\chi(z)^2}{\tau} \bigg{(} 
  \Phi_j (\alpha, \omega(z)) - 
  \Phi_j(\alpha-V, \omega(z)) \bigg{)}.
$$
Note that 
\begin{eqnarray} \label{eq_rho}
 |\rho_j (\tau,z)| 
 \,\le\, 
 C\eps^3\tau^{-2} 
 +C\eps \tau^{-2+3\delta/2} 
 +C\eps^5 \tau^{-2} 
 \,\le\, 
 C\eps\tau^{-2+3\delta/2}
 \,\le\,
 C\eps \tau^{-\eta},
\end{eqnarray}
where $2>\eta=2-3\delta/2>1$, since each term of
$$
 \Fjcub (\realpart(\alpha e^{i \theta}), \cdots) 
 -\Fjcub (\realpart((\alpha-V) e^{i \theta}), \cdots) 
$$
includes at least one or more $V_j$ 
which has better time decay as compared to $\alpha_j$.
We also note that the condition \eqref{eq_decayassumption1} implies
\begin{eqnarray} \label{eq_Phicondi}
 \imagpart \jb{\Phi(\beta(\tau,z), \omega(z)), A\beta(\tau,z)}_{\C^N} 
 \le 
 - C^* \omega_0(z) |\beta(\tau,z)|^4
\end{eqnarray}
for some positive constant $C^*$.
Now, similarly to \cite{KLS}, we put 
$\Psi(\tau)=\nu_A(\beta(\tau,z))^2$,
$\rho = (\rho_j)_{1\le j\le N}$
and compute
$$
 \frac{d}{d\tau} \left( (\log (\tau \omega_0(z)))^2 \Psi(\tau) \right) 
 = 
 (\log (\tau \omega_0(z)))^2 \frac{d \Psi}{d\tau}(\tau) 
 +\frac{2\log(\tau \omega_0(z))}{\tau} \Psi(\tau).
$$
By \eqref{eq_betaDE}, \eqref{eq_rho} and \eqref{eq_Phicondi}, we have
\begin{eqnarray*}
 \frac{d \Psi}{d \tau} (\tau) 
 &=& 
 2 \realpart \jb{\pa_\tau \beta(\tau,z), A\beta(\tau,z)}_{\C^N} \\
 &=& 
 \frac{2 \chi(z)^2}{\tau} 
  \imagpart \jb{\Phi(\beta(\tau,z), \omega(z)), A\beta(\tau,z)}_{\C^N} 
  +2\realpart \jb{\rho(\tau,z), A \beta(\tau,z)}_{\C^N} \\
 &\le& 
 -\frac{2}{\tau}C^* \chi(z)^2 \omega_0(z) |\beta(\tau,z)|^4 
 +\frac{C \eps^2}{\tau^\eta},
\end{eqnarray*}
where $C^*$ is the constant appearing in \eqref{eq_Phicondi}. 
We also have 
\begin{eqnarray*}
 \Psi(\tau) 
 &\le&  \lambda^*  |\beta(\tau,z)|^2 \\
 &=& 
 \frac{\lambda^*}
      {\sqrt{2C^*\chi(z)^2\omega_0(z)\log(\tau \omega_0(z))}} 
   \sqrt{2C^*\chi(z)^2 \omega_0(z)\log(\tau\omega_0(z))} |\beta(\tau,z)|^2 
  \\
 &\le& 
 \frac{(\lambda^*)^2}{4C^* \chi(z)^2\omega_0(z)\log(\tau\omega_0(z))} 
 +
 C^*\chi(z)^2\omega_0(z)\log(\tau\omega_0(z)) |\beta(\tau,z)|^4.
\end{eqnarray*}
Piecing them together, we obtain
$$
 \frac{d}{d \tau}\left((\log(\tau \omega_0(z)))^2 \Psi(\tau) \right) 
 \le 
 \frac{C}{\tau \chi(z)^2 \omega_0(z)} 
 +\frac{C\eps^2(\log(\tau \omega_0(z)))^2}{\tau^\eta}.
$$
Integrating with respect to $\tau$, we arrive at
\begin{align*}
 &(\log(\tau \omega_0(z)))^2 \Psi(\tau) \\
 &\le
 (\log(\tau_0 \omega_0(z)))^2 \nu_A(\beta(\tau_0,z))^2
 +\int_{\tau_0}^\tau
  \frac{C}{\sigma\chi(z)^2\omega_0(z)}
  +\frac{C\eps^2(\log(\sigma\omega_0(z)))^2}{\sigma^\eta}\,d\sigma\\
 &\le
 (\log(\tau_0 \omega_0(z)))^2 C \eps^2
  +\frac{C\log(\tau\omega_0(z))}{\chi(z)^2\omega_0(z)}
  +C\eps^2\int_1^\infty \omega_0(z)^{\eta-1} \frac{(\log\sigma)^2}{\sigma^\eta}
 \,d\sigma\\
 &\le
 \frac{C\log(\tau \omega_0(z))}{\chi(z)^2 \omega_0(z)} 
 +C\eps^2 \omega_0(z),
\end{align*}
whence
$$
 |\beta(\tau,z)| 
 \le 
 \frac{1}{\sqrt{\lambda_*}} \Psi(\tau)^{1/2} 
 \le 
 \frac{C}{\sqrt{\chi(z)^2 \omega_0(z) \log(\tau \omega_0(z))}}
 +\frac{C\eps \omega_0(z)^{1/2}}{\log(\tau \omega_0(z))}
$$
for $\tau \ge \tau_0$. Therefore 
we have
\begin{align} 
 |\alpha(\tau,z)| \chi(z) \omega_0(z)^{1/2}
 \le (|\beta(\tau,z)|+|V(\tau,z)|)\chi(z) \omega_0(z)^{1/2}
 \le 
 \frac{C}{\sqrt{\log(\tau \omega_0(z))}}.
 \label{eq_inparticular}
\end{align}
It follows from  \eqref{eq_substitution} and \eqref{eq_inparticular} that
\begin{eqnarray} \label{eq_udecay}
 |u(t,x)| 
 \,\le\, 
 \frac{|\alpha(\tau,z)| \chi(z) \omega_0(z)^{1/2}}
       {\sqrt{\tau \omega_0(z)}} 
 \,\le\,
 \frac{C}{\sqrt{\tau \omega_0(z) \log(\tau \omega_0(z))}} 
 \,\le\,
 \frac{C(1+t)^{-1/2}}{\sqrt{\log(2+t)}}.
\end{eqnarray}
Using \eqref{eq_udecay} and the finite propagation speed, we obtain
\begin{eqnarray*}
 \|u(t, \cdot) \|_{L^p(\R)} 
 &\le&
 C\frac{(1+t)^{-(1/2-1/p)}}{\sqrt{\log(2+t)}}
\end{eqnarray*}
for $p \in [2,\infty]$, 
which proves Theorem~\ref{thm_decay}.\\

To prove Theorem~\ref{thm_decay2}, 
we remark that the condition \eqref{eq_decayassumption2} implies
\begin{eqnarray*}
 \imagpart \jb{\Phi(\beta(\tau,z), \omega(z)), A\beta(\tau,z)}_{\C^N} 
 \le 
 - C_* \omega_0(z)^3 |\beta(\tau,z)|^4
\end{eqnarray*}
for some positive constant $C_*$. Repeating the same steps as before,
we get
$$
 |\alpha(\tau,z)| \chi(z) \omega_0(z)^{3/2}
 \le 
 \frac{C}{\sqrt{\log(\tau \omega_0(z))}}
$$
instead of \eqref{eq_inparticular}.
Finally, we obtain
\begin{eqnarray*}
 \sum_{|I|=1} \left| \pa_{t,x}^I u(t,x) \right|
 &\le&
 \frac{C}{\sqrt{\tau\omega_0(z) \log(\tau \omega_0(z))}}
 +\frac{C\chi(z)\omega_0(z)^{1+4/3}}{(\tau\omega_0(z))^{4/3}}
 +\frac{C\chi(z)\omega_0(z)^{1+3/2}}{(\tau\omega_0(z))^{3/2}}\\
 &\le&
 \frac{C(1+t)^{-1/2}}{\sqrt{\log(2+t)}}
\end{eqnarray*}
and
$$
 \sum_{|I|\le 1} \| \pa_{t,x}^I u(t, \cdot) \|_{L^p(\R)}
 \le
 \frac{C(1+t)^{-(1/2-1/p)}}{\sqrt{\log(2+t)}}
$$
for $p \in [2,\infty ]$, 
which completes the proof of Theorem~\ref{thm_decay2}.
\qed\\

\medskip
\subsection*{Acknowledgment}
The work of H.~S. 
is supported by Grant-in-Aid for Scientific Research (C) (No.~25400161), JSPS.


\end{document}